\documentclass{amsart}
\usepackage{amsmath,amssymb,amsfonts}
\usepackage[mathscr]{eucal}

\usepackage{enumerate}
\newenvironment{enumroman}{\begin{enumerate}[\upshape (i)]}
                                                {\end{enumerate}}

\theoremstyle{plain}
\newtheorem{theorem}{Theorem}[section]

\newtheorem{prop}[theorem]{Proposition}
\newtheorem{lemma}[theorem]{Lemma}

\theoremstyle{definition}
\newtheorem{definition}[theorem]{Definition}

\newtheorem*{thank}{Acknowledgments}

\input xy
\xyoption{all}

\newcommand{\Deltaop}{{\bf \Delta}^{op}}
\newcommand{\IDeltaop}{{\bf I \Delta}^{op}}
\newcommand{\colim}{\text{colim}}
\newcommand{\hocolim}{\text{hocolim}}

\newcommand{\nerve}{\text{nerve}}
\newcommand{\Hom}{\text{Hom}}
\newcommand{\Map}{\text{Map}}
\newcommand{\SSets}{\mathcal{SS}ets}
\newcommand{\SSetsd}{\mathcal{SS}ets^{\mathcal D}}
\newcommand{\SSetst}{\mathcal{SS}ets^{\mathcal T}}
\newcommand{\SSetstm}{\mathcal{SS}ets^{\mathcal T_M}}
\newcommand{\LSSetst}{\mathcal {LSS}ets^{\mathcal T}}
\newcommand{\LSSetstm}{\mathcal {LSS}ets^{\mathcal T_M}}
\newcommand{\Tm}{\mathcal T_M}
\newcommand{\Sets}{\mathcal Sets}
\newcommand{\Algt}{\mathcal Alg^\mathcal T}
\newcommand{\Algtm}{\mathcal Alg^{\mathcal T_M}}
\newcommand{\Algtam}{\mathcal Alg^{\mathcal T_{AM}}}
\newcommand{\SSpo}{\mathcal{SS} p_\mathcal O}
\newcommand{\SSpof}{\mathcal{SS} p_{\mathcal O,f}}
\newcommand{\SSpoc}{\mathcal{SS} p_{\mathcal O,c}}
\newcommand{\LSSpof}{\mathcal{LSS}p_{\mathcal O,f}}
\newcommand{\LSSpoc}{\mathcal{LSS}p_{\mathcal O,c}}
\newcommand{\SSp}{\mathcal{SS}p}
\newcommand{\LSSp}{\mathcal{LSS}p}
\newcommand{\alphau}{{\underline \alpha}}
\newcommand{\xu}{{\underline x}}
\newcommand{\Tocat}{\mathcal T_{\mathcal {OC}at}}
\newcommand{\SSetsTocat}{\mathcal {SS}ets^{\mathcal
T_{\mathcal {OC}at}}}
\newcommand{\LSSetsTocat}{\mathcal {LSS}ets^{\mathcal
T_{\mathcal {OC}at}}}

\newcommand{\ISSpof}{\mathcal{ISS}p_{\mathcal O,f}}
\newcommand{\ISSpoc}{\mathcal{ISS}p_{\mathcal O,c}}
\newcommand{\ILSSpof}{\mathcal{ILSS}p_{\mathcal O,f}}
\newcommand{\ILSSpoc}{\mathcal{ILSS}p_{\mathcal O,c}}
\newcommand{\SSetstg}{\mathcal {SS}ets^{\mathcal T_G}}
\newcommand{\LSSetstg}{\mathcal {LSS}ets^{\mathcal T_G}}
\newcommand{\Tg}{\mathcal T_G}
\newcommand{\Algtg}{\mathcal Alg^{\mathcal T_G}}

\newcommand{\ILSSp}{\mathcal{ILSS}p}
\newcommand{\Togd}{\mathcal T_{\mathcal {OG}d}}

\newcommand{\LSSetsTogd}{\mathcal {LSS}ets^{\mathcal
T_{\mathcal {OG}d}}}
\newcommand{\Gammaop}{\Gamma^{op}}
\newcommand{\Tam}{\mathcal T_{AM}}
\newcommand{\Algtgd}{\mathcal Alg^{\mathcal T_{Gd}}}
\newcommand{\Algtogd}{\mathcal Alg^{\mathcal T_{\mathcal O Gd}}}
\newcommand{\Algtocat}{\mathcal Alg^{\mathcal T_{\mathcal O Cat}}}

\begin{document}

\title[Adding Inverses]{Adding Inverses to Diagrams Encoding Algebraic Structures}

\author{Julia E. Bergner}

\email{bergnerj@member.ams.org}

\address{Department of Mathematics, Kansas State University, 138 Cardwell
Hall, Manhattan, KS 66506}

\curraddr{University of California, Riverside, CA 92521}

\subjclass[2000]{55U10, 18B40, 18C10, 55P35}

\keywords{simplicial groups, Segal groupoids, diagram categories}

\begin{abstract}
We modify a previous result, which showed that certain diagrams of
spaces are essentially simplicial monoids, to construct diagrams
of spaces which model simplicial groups.  Furthermore, we show
that these diagrams can be generalized to models for Segal
groupoids. We then modify Segal's model for simplicial abelian
monoids in such a way that it becomes a model for simplicial
abelian groups.
\end{abstract}

\maketitle

\section{Introduction}

Much research has been done on various structures equivalent to
topological or simplicial groups.  Classical examples include
Thomason's work with delooping machines \cite{th} and Stasheff's
group-like $A_\infty$-spaces \cite{stash}, while more modern
examples make use of the structure of algebraic theories
\cite{bad}, \cite{bcv}.  In \cite{simpmon}, we focus instead on
simplicial monoids.  In particular, we consider diagrams of spaces
which essentially give one of the spaces the structure of a
monoid.  While such a structure is given when the diagram is the
theory of monoids, we obtain a simpler diagram which, from the
viewpoint of homotopy theory, encodes the same structure.

Here we return to the original situation and ask if this
construction can be modified to obtain a diagram encoding the
structure of a simplicial group rather than a simplicial monoid.
After all, we need only to find a way to include inverses.  In
this paper, we show that indeed we can represent simplicial groups
in this manner.

Furthermore, as the construction for simplicial monoids
generalizes to the many-object case of simplicial categories, we
can obtain models for simplicial groupoids as well.  In doing so,
we extend a result relating simplicial categories to Segal
categories (essentially categories with composition only given up
to higher homotopy) to one relating simplicial groupoids to Segal
groupoids.

One might ask, then, if there are other algebraic structures which
can be modelled in a similar way, in particular by some kind of
diagram which is simpler than the one given by the corresponding
algebraic theory.  An answer for the case of simplicial abelian
monoids is given by Segal's category $\Gamma$ \cite{segal},
although it does not provide the relationship between strict and
homotopy structures that we have in the non-abelian situation.
Here, we extend this construction to find a model for abelian
groups.

To obtain this result, we use work of Bousfield, in which he
modifies the projections in the diagram for simplicial monoids,
rather than the diagram itself, to encode inverses and so model
simplicial groups \cite{bous}. To model simplicial abelian groups,
we change the projections in the category $\Gamma$ encoding an
abelian monoid structure in an analogous way.  It should be noted
that Bousfield's work includes models for $n$-fold loop spaces,
work which has been approached from another angle by Berger
\cite{berger}.

Of course, the question we answer here in these cases can be asked
for many other algebraic structures.  For example, are there
simple diagrams encoding the structure of (commutative) rings?  We
can also ask about so-called multi-sorted algebraic structures,
such as operads and group actions, for which we have models given
by multi-sorted algebraic theories \cite{multisort}. We hope to
find more such results in future work.  In related work, Barwick
has found diagrams given by operator categories which model
certain kinds of algebraic structures \cite{bar}.

In sections 2 and 3, we review basic properties of model
categories and some of the model category structures that we use
in this paper.  In section 4, we modify the argument given in
\cite{simpmon} for simplicial monoids to obtain a model for
simplicial groups. We extend this argument in section 5 to
simplicial groupoids with a fixed object set.  In section 6, we
summarize the construction of Bousfield's alternative model for
simplicial groups.  Finally, in section 7, we summarize Segal's
argument for a model for simplicial abelian monoids and modify it
to obtain a model for simplicial abelian groups.

\begin{thank}
I would like to thank Bernard Badzioch and Bill Dwyer for their
helpful suggestions for this paper, particularly with regard to
the approach to simplicial abelian groups, and Clark Barwick and
Joachim Kock for conversations about the category $\IDeltaop$. I
am also grateful to the referee for helpful comments and
corrections.
\end{thank}

\section{Model Category Structures}

In this section, we review the necessary tools from model
categories that we will need to prove our result for simplicial
groups.

Recall that a model category structure on a category $\mathcal C$
is a choice of three distinguished classes of morphisms:
fibrations, cofibrations, and weak equivalences.  A (co)fibration
which is also a weak equivalence will be called an \emph{acyclic
(co)fibration}.  With this choice of three classes of morphisms,
$\mathcal C$ is required to satisfy five axioms MC1-MC5
\cite[3.3]{ds}.  Here we only state MC4 and MC5, as they are used
in the course of the paper.

\begin{itemize}
%
%
%
\item (MC4) If $i:A \rightarrow B$ is a cofibration and $p:X
\rightarrow Y$ is a fibration, then a dotted arrow lift exists in
any solid arrow diagram of the form
\[ \xymatrix{A \ar[r] \ar[d]^i & X \ar[d]^p \\
B \ar[r] \ar@{-->}[ur] & Y} \] if either
\begin{enumroman}
\item $p$ is a weak equivalence, or

\item $i$ is a weak equivalence.
\end{enumroman}
(In this case we say that $i$ has the \emph{left lifting property}
with respect to $p$ and that $p$ has the \emph{right lifting
property} with respect to $i$.)

\item (MC5) Any map $f$ can be factored two ways:
\begin{enumroman}
\item $f= pi$ where $i$ is a cofibration and $p$ is an acyclic
cofibration, and

\item $f=qj$ where $j$ is an acyclic cofibration and $q$ is a
fibration.
\end{enumroman}
\end{itemize}

An object $X$ in $\mathcal C$ is \emph{fibrant} if the unique map
$X \rightarrow \ast$ from $X$ to the terminal object is a
fibration.  Dually, $X$ is \emph{cofibrant} if the unique map
$\phi \rightarrow X$ from the initial object to $X$ is a
cofibration.  The factorization axiom MC5 guarantees that each
object $X$ has a weakly equivalent fibrant replacement $\widehat
X$ and a weakly equivalent cofibrant replacement $\widetilde X$.
These replacements are not necessarily unique, but they can be
chosen to be functorial in the cases we will use
\cite[1.1.3]{hovey}.

The model category structures which we will discuss are all
cofibrantly generated.  In a cofibrantly generated model category,
there are two sets of morphisms, one of generating cofibrations
and one of generating acyclic cofibrations, such that a map is a
fibration if and only if it has the right lifting property with
respect to the generating acyclic cofibrations, and a map is an
acyclic fibration if and only if it has the right lifting property
with respect to the generating cofibrations \cite[11.1.2]{hirsch}.
The following theorem is useful for proving the existence of a
cofibrantly generated model structure on a category.

\begin{theorem} \cite[11.3.1]{hirsch} \label{CofGen}
Let $\mathcal M$ be a category which has all small limits and
colimits. Suppose that $\mathcal M$ has a class of weak
equivalences which satisfies the two-out-of-three property (model
category axiom MC2) and which is closed under retracts. Let $I$
and $J$ be sets of maps in $\mathcal M$ which satisfy the
following conditions:
\begin{enumerate}
\item Both $I$ and $J$ permit the small object argument
\cite[10.5.15]{hirsch}.

\item Every $J$-cofibration is an $I$-cofibration and a weak
equivalence.

\item Every $I$-injective is a $J$-injective and a weak
equivalence.

\item One of the following conditions holds:
\begin{enumerate}
\item A map that is an $I$-cofibration and a weak equivalence is a
$J$-cofibration, or

\item A map that is both a $J$-injective and a weak equivalence is
an $I$-injective.
\end{enumerate}
\end{enumerate}
Then there is a cofibrantly generated model category structure on
$\mathcal M$ in which $I$ is a set of generating cofibrations and
$J$ is a set of generating acyclic cofibrations.
\end{theorem}

We now state the definition of a Quillen pair of model category
structures.  Recall that for categories $\mathcal C$ and $\mathcal
D$ a pair of functors
\[ \xymatrix@1{F: \mathcal C \ar@<.5ex>[r] & \mathcal D:R
\ar@<.5ex>[l]} \] is an \emph{adjoint pair} if for each object $X$
of $\mathcal C$ and object $Y$ of $\mathcal D$ there is a natural
isomorphism $\psi :\Hom_\mathcal D(FX,Y) \rightarrow \Hom_\mathcal
C(X,RY)$ \cite[IV.1]{macl}.  The adjoint pair is sometimes written
as the triple $(F, R, \psi)$.

\begin{definition} \cite[1.3.1]{hovey}
If $\mathcal C$ and $\mathcal D$ are model categories, then an
adjoint pair $(F,R, \psi)$ between them is a \emph{Quillen pair}
if one of the following equivalent statements holds:
\begin{enumerate}
\item $F$ preserves cofibrations and acyclic cofibrations.

\item $R$ preserves fibrations and acyclic fibrations.
\end{enumerate}
\end{definition}

We now have the following definition of Quillen equivalence, which
is the standard notion of equivalence of model category
structures.

\begin{definition} \cite[1.3.12]{hovey}
A Quillen pair is a \emph{Quillen equivalence} if for all
cofibrant $X$ in $\mathcal C$ and fibrant $Y$ in $\mathcal D$, a
map $f:FX \rightarrow Y$ is a weak equivalence in $\mathcal D$ if
and only if the map $\psi f:X \rightarrow RY$ is a weak
equivalence in $\mathcal C$.
\end{definition}

We will use the following proposition to prove that our Quillen
pairs are Quillen equivalences.  Recall that a functor $F:\mathcal
C \rightarrow \mathcal D$ \emph{reflects} a property if, for any
morphism $f$ of $\mathcal C$, whenever $Ff$ has the property, then
so does $f$.

\begin{prop}\cite[1.3.16]{hovey}
Suppose that $(F,R,\psi)$ is a Quillen pair from $\mathcal C$ to
$\mathcal D$. Then the following statements are equivalent:
\begin{enumerate}
\item $(F,R, \psi)$ is a Quillen equivalence.


\item $R$ reflects weak equivalences between fibrant objects, and
for every cofibrant $X$ in $\mathcal C$ the map $X \rightarrow
R(FX)^f$ is a weak equivalence.
\end{enumerate}
\end{prop}

Throughout this paper, we use the category of simplicial sets,
denoted $\SSets$.  Recall that a simplicial set is a functor
$\Deltaop \rightarrow \Sets$, where ${\bf \Delta}$ denotes the
cosimplicial category whose objects are the finite ordered sets
$[n]=(0 \rightarrow \cdots \rightarrow n)$ and whose morphisms are
the order-preserving maps.  The simplicial category $\Deltaop$ is
then the opposite of this category.  Some examples of simplicial
sets are, for each $n \geq 0$, the $n$-simplex $\Delta [n]$, its
boundary $\dot \Delta [n]$, and, for any $0 \leq k \leq n$, the
simplicial set $V[n,k]$, which is $\dot \Delta [n]$ with the $k$th
face removed \cite[I.1]{gj}.  More generally, a \emph{simplicial
object} in a category $\mathcal C$ is a functor $\Deltaop
\rightarrow \mathcal C$.  In particular, a functor $\Deltaop
\rightarrow \SSets$ is a \emph{simplicial space} or
\emph{bisimplicial set} \cite[IV]{gj}.

In a slight abuse of terminology, we use the term \emph{simplicial
category} to refer to a simplicial object in the category of all
(small) categories which satisfies the additional condition that
the face and degeneracy maps are the identity on all the objects.
Such an object is often called a \emph{category enriched over
simplicial sets}, since it is just a category with a simplicial
set of morphisms between any two objects.  For any objects $X$ and
$Y$ in a simplicial category, we denote this simplicial set
$\Map(X,Y)$ and call it a \emph{function complex}.  A simplicial
category in this sense with all morphisms invertible is called a
\emph{simplicial groupoid}.

%

We also use the notion of a simplicial model category $\mathcal
M$, or a model category which is also a simplicial category
satisfying two axioms \cite[9.1.6]{hirsch}. It is important to
note that a function complex in a simplicial model category is
only homotopy invariant in the case that $X$ is cofibrant and $Y$
is fibrant. For the general case, we have the following
definition:

\begin{definition} \cite[17.3.1]{hirsch}
A \emph{homotopy function complex} $\Map^h(X,Y)$ in a simplicial
model category $\mathcal M$ is the simplicial set $\Map(\widetilde
X, \widehat Y)$ where $\widetilde X$ is a cofibrant replacement of
$X$ in $\mathcal M$ and $\widehat Y$ is a fibrant replacement for
$Y$.
\end{definition}

Several of the model category structures that we use are obtained
by localizing a given model category structure with respect to a
map or a set of maps.  Suppose that $P = \{f:A \rightarrow B\}$ is
a set of maps with respect to which we would like to localize a
model category $\mathcal M$.

\begin{definition} \label{local}
A $P$-\emph{local object} $W$ is a fibrant object of $\mathcal M$
such that for any $f:A \rightarrow B$ in $P$, the induced map on
homotopy function complexes
\[ f^*:\Map^h(B,W) \rightarrow \Map^h(A,W) \]
is a weak equivalence of simplicial sets.  A map $g:X \rightarrow
Y$ in $\mathcal M$ is then a $P$-\emph{local equivalence} if for
every local object $W$, the induced map on homotopy function
complexes
\[ g^*: \Map^h(Y,W) \rightarrow \Map^h(X,W) \]
is a weak equivalence of simplicial sets.
\end{definition}

\begin{theorem}\cite[4.1.1]{hirsch} \label{localize}
Let $\mathcal M$ be a left proper cellular model category
\cite[13.1.1, 12.1.1]{hirsch} and $P$ a set of morphisms of
$\mathcal M$. There is a model category structure $\mathcal L_P
\mathcal M$ on the underlying category of $\mathcal M$ such that:
\begin{enumerate}
\item The weak equivalences are the $P$-local equivalences.

\item The cofibrations are precisely the cofibrations of $\mathcal
M$.

\item The fibrations are the maps which have the right lifting
property with respect to the maps which are both cofibrations and
$P$-local equivalences.

\item The fibrant objects are the $P$-local objects.
\end{enumerate}
\end{theorem}

In this situation, we refer to the functorial fibrant replacement
functor as a \emph{localization} functor.

When we are working with localized model category structures, the
following theorem can be used to prove that an adjoint pair is
still a Quillen pair after localization.

\begin{theorem} \cite[3.3.20]{hirsch} \label{LocPair}
Let $\mathcal C$ and $\mathcal D$ be left proper, cellular model
categories and let $(F,R, \psi)$ be a Quillen pair between them.
Let $S$ be a set of maps in $\mathcal C$ and $L_S \mathcal C$ the
localization of $\mathcal C$ with respect to $S$.  Then if ${\bf
L}FS$ is the set in $\mathcal D$ obtained by applying the left
derived functor of $F$ to the set $S$ \cite[8.5.11]{hirsch}, then
$(F,R, \psi)$ is also a Quillen pair between the model categories
$L_S\mathcal C$ and $L_{{\bf L}FS}D$.
\end{theorem}

One important model category structure we will use is the standard
one on the category $\SSets$ of simplicial sets.  In this case, a
weak equivalence is a map of simplicial sets $f:X \rightarrow Y$
such that the induced map $|f|:|X| \rightarrow |Y|$ is a weak
homotopy equivalence of topological spaces.  The cofibrations are
monomorphisms, and the fibrations are the maps with the right
lifting property with respect to the acyclic cofibrations
\cite[I.11.3]{gj}. This model category structure is cofibrantly
generated; a set of generating cofibrations is $I=\{\dot \Delta
[n] \rightarrow \Delta [n] \mid n \geq 0\}$, and a set of
generating acyclic cofibrations is $J=\{V[n,k] \rightarrow \Delta
[n] \mid n \geq 1, 0 \leq k \leq n\}$.

\section{Model categories of diagrams}

The objects of all the categories we consider in this paper are
given by diagrams of spaces.  Given any small category $\mathcal
D$, there is a category $\SSetsd$ of $\mathcal D$-diagrams in
$\SSets$, or functors $\mathcal D \rightarrow \SSets$.  We can
obtain two model category structures on $\SSetsd$ by the following
results.

\begin{theorem} \cite[IX 1.4]{gj}
Given the category $\mathcal{SS}ets^\mathcal D$ of $\mathcal
D$-diagrams of simplicial sets, there is a simplicial model
category structure $\mathcal{SS}ets^\mathcal D_f$ in which the
weak equivalences and fibrations are objectwise and in which the
cofibrations are the maps which have the left lifting property
with respect to the maps which are both fibrations and weak
equivalences.
\end{theorem}

\begin{theorem} \cite[VIII 2.4]{gj}
There is a simplicial model category $\mathcal{SS}ets^\mathcal
D_c$ in which the weak equivalences and the cofibrations are
objectwise and in which the fibrations are the maps which have the
right lifting property with respect to the maps which are
cofibrations and weak equivalences.
\end{theorem}

Given these general results, we now turn to the particular
diagrams which we will be considering in this paper.  Let $G$ be a
simplicial group, by which we mean a simplicial object in the
category $\mathcal Grp$ of groups, or a functor $\Deltaop
\rightarrow \mathcal Grp$.  However, here we use an alternate
viewpoint in which we use algebraic theories to define simplicial
groups. We begin with the definition of an algebraic theory. Some
references for algebraic theories include chapter 3 of \cite{bor},
the introduction to \cite{bad}, and section 3 of \cite{multisort}.

\begin{definition}
An \emph{algebraic theory} $\mathcal T$ is a small category with
finite products and objects denoted $T_n$ for $n \geq 0$. For each
$n$, $T_n$ is equipped with an isomorphism $T_n \cong (T_1)^n$.
Note in particular that $T_0$ is the terminal object in $T$.
\end{definition}

Here we consider one particular theory, the theory of groups,
which we denote $\Tg$. To describe this theory, we first consider
the full subcategory of the category of groups generated by
representatives $T_n$ of the isomorphism classes of free groups on
$n$ generators. We then define the \emph{theory of groups} $\Tg$
to be the opposite of this category. Thus $T_n$, which is
canonically the coproduct of $n$ copies of $T_1$ in the category
$\mathcal Grp$ of groups, becomes the product of $n$ copies of
$T_1$ in $\Tg$.  It follows that there is a projection map
$p_{n,i}: T_n \rightarrow T_1$ for each $1 \leq i \leq n$ in
addition to other group homomorphisms. In fact, there are such
projection maps in any algebraic theory. We use them to make the
following definition.

\begin{definition} \cite[1.1]{bad}
Given an algebraic theory $\mathcal T$, a \emph{strict simplicial}
$\mathcal T$-\emph{algebra} (or simply $\mathcal
T$-\emph{algebra}) $A$ is a product-preserving functor $A:\mathcal
T \rightarrow \mathcal {SS}ets$.  Here, ``product-preserving"
means that for each $n \geq 0$ the canonical map
\[ A(T_n) \rightarrow A(T_1)^n, \]
induced by the $n$ projection maps $T_n \rightarrow T_1$, is an
isomorphism of simplicial sets.  In particular, $A(T_0)$ is the
one-point simplicial set $\Delta [0]$.
\end{definition}

In general, a $\mathcal T$-algebra $A$ defines a strict algebraic
structure on the space $A(T_1)$ corresponding to the theory
$\mathcal T$ \cite[\S 1]{bad}. So, a $\Tg$-algebra $A$ defines a
group structure on the space $A(T_1)$.  In fact, the category of
simplicial groups is equivalent to the category of $\Tg$-algebras
\cite{law}, \cite[II.4]{quillen}.

We can also consider the case where the products are not preserved
strictly, but only up to homotopy.

\begin{definition} \cite[1.2]{bad}
Given an algebraic theory $\mathcal T$, a \emph{homotopy}
$\mathcal T$-\emph{algebra} is a functor $X:\mathcal T \rightarrow
\SSets$ which preserves products up to homotopy.  The functor $X$
preserves products up to homotopy if, for each $n \geq 0$ the
canonical map
\[ X(T_n) \rightarrow X(T_1)^n \]
induced by the projection maps $p_{n,i}:T_n \rightarrow T_1$ for
$1 \leq i \leq n$ is a weak equivalence of simplicial sets. In
particular, we assume that $X(T_0)$ is weakly equivalent to
$\Delta [0]$.
\end{definition}

We now consider the corresponding model category structures.

\begin{prop} \cite[II.4]{quillen}, \cite[3.1]{sch}
Let $\mathcal T$ be an algebraic theory and $\Algt$ the category
of $\mathcal T$-algebras.  Then there is a cofibrantly generated
model category structure on $\Algt$ in which the weak equivalences
and fibrations are levelwise weak equivalences of simplicial sets
and the cofibrations are the maps with the left lifting property
with respect to the maps which are fibrations and weak
equivalences.
\end{prop}

We also need a model category structure for homotopy $\mathcal
T$-algebras. However, there is no model category structure on the
category of homotopy $\mathcal T$-algebras itself since this
category does not have all small colimits. However, there is a
model category structure on the category of all $\mathcal
T$-diagrams of simplicial sets in which the fibrant objects are
homotopy $\mathcal T$-algebras.  To obtain this structure, we
begin by considering the model category structure $\SSetst_f$ on
the category of all functors $\mathcal T \rightarrow \SSets$. The
desired model structure can be obtained by localizing the model
structure $\SSetst_f$ with respect to a set of maps. We summarize
this localization here; a complete description is given by
Badzioch \cite[\S 5]{bad}.

Given an algebraic theory $\mathcal T$, consider the functor
$\Hom_\mathcal T (T_n,-)$.  We then have maps
\[ p_n: \coprod_{i=1}^n \Hom_\mathcal T(T_1,-) \rightarrow \Hom_\mathcal T(T_n,-) \]
induced from the projection maps in $\mathcal T$.  We then
localize the model category structure on $\SSetst$ with respect to
the set $S= \{p_n \mid n \geq 0 \}$. We denote the resulting model
category structure $\LSSetst$.

\begin{prop} \cite[5.5]{bad}.
The fibrant objects in $\LSSetst$ are the homotopy $\mathcal
T$-algebras which are fibrant in $\SSetst$.
\end{prop}

We now have the following result by Badzioch which relates strict
and homotopy $\mathcal T$-algebras.

\begin{theorem}  \cite[6.4]{bad} \label{rigid}
Given an algebraic theory $\mathcal T$, there is a Quillen
equivalence of model categories between $\Algt$ and $\LSSetst$.
\end{theorem}

We will find it convenient, however, to work in the situation
where a homotopy $\mathcal T$-algebra $X$ has $X_0$ precisely
$\Delta[0]$ rather than just a space weakly equivalent to it.  The
following two results are proved in \cite{simpmon} and
\cite{simpmon2} for the theory of monoids $\Tm$, but their proofs
hold for any algebraic theory $\mathcal T$.

\begin{prop} \cite[3.11]{simpmon}
Consider the category $\SSetst_*$ of functors $\mathcal T
\rightarrow \SSets$ such that the image of $T_0$ is $\Delta [0]$.
There is a model category structure on $\SSetst_*$ in which the
weak equivalences and fibrations are defined levelwise and the
cofibrations are the maps with the left lifting property with
respect to the acyclic fibrations.
\end{prop}

Now, to obtain a localized model category $\LSSetst_*$, we need to
modify the maps
\[ p_n: \coprod_{i=1}^n \Hom_\mathcal T(T_1,-) \rightarrow \Hom_\mathcal
T(T_n,-) \] that we used to obtain $\LSSetst$ from $\SSetst$.
Since $\Hom_\mathcal T(T_n,-)_0 \cong \Delta [0]$ for all $n \geq
0$, the only change we need to make to these maps is to take the
coproduct $\coprod_n \Hom_\mathcal T(T_1,-)$ in the category
$\SSetstm_*$ (as in \cite[3.6]{simpmon}).  We then localize
$\SSetstm_*$ with respect to the set of all such maps to obtain a
model structure $\LSSetst_*$.

Since a fibrant and cofibrant object $X$ in $\LSSetst$ has $X_0$
weakly equivalent to $\Delta [0]$, it is not too surprising that
we can instead use $\LSSetst_*$ in the following variation on
Badzioch's result:

\begin{prop} \cite{simpmon2}
There is a Quillen equivalence of model categories
\[ \xymatrix@1{L:\LSSetst_* \ar@<.5ex>[r] & \Algt:I, \ar@<.5ex>[l]}
\] where $I$ denotes the inclusion functor and $L$ is its left
adjoint.
\end{prop}

Now that we have established our definitions for simplicial
groups, we turn to reduced Segal categories and reduced Segal
groupoids. These are simplicial spaces satisfying some conditions,
so, like the simplicial monoids, they are given by a diagram of
simplicial sets. We begin with the definition of a Segal
precategory.

\begin{definition}
A \emph{Segal precategory} is a simplicial space $X$ such that
$X_0$ is a discrete simplicial set.  If $X_0$ consists of a single
point, then $X$ is a \emph{reduced Segal precategory}.
\end{definition}

Now note that for any simplicial space $X$ there is a \emph{Segal
map}
\[ \varphi_n: X_n \rightarrow \underbrace{X_1 \times_{X_0} \ldots
\times_{X_0} X_1}_n \] for each $n \geq 2$, which we define as
follows.  Let $\alpha^k:[1] \rightarrow [n]$ be the map in ${\bf
\Delta}$ such that $\alpha^k(0)=k$ and $\alpha^k(1)=k+1$, defined
for each $0 \leq k \leq n-1$. We can then define the dual maps
$\alpha_k:[n]\rightarrow [1]$ in $\Deltaop$.  For $n \geq 2$, the
Segal map is defined to be the map
\[ \varphi_n: X_n \rightarrow \underbrace{X_1 \times_{X_0} \cdots \times_{X_0} X_1}_n \]
induced by the maps
\[ X(\alpha_k):X_n \rightarrow X_1. \]

\begin{definition} \cite[\S 2]{hs} \label{SeCat}
A \emph{Segal category} $X$ is a Segal precategory $X: \Deltaop
\rightarrow \SSets$ such that $X_0= \Delta [0]$ and such that for
each $n \geq 2$ the Segal map
\[ \varphi_n: X_n \rightarrow \underbrace{X_1 \times_{X_0} \cdots \times_{X_0} X_1}_n \]
is a weak equivalence of simplicial sets.  If $X_0= \Delta [0]$,
then $X$ is a \emph{reduced Segal category}.  Note that in this
case we have $X_n \simeq (X_1)^n$.
\end{definition}

We now give model category structures on the category $\SSpo$ of
Segal precategories with a fixed object set.  As in the case of
homotopy $\mathcal T$-algebras, there is no model structure on the
category of Segal categories, due to a lack of colimits.

\begin{prop} \cite[3.7]{simpmon}
There is a model category structure on $\SSpo$, which we denote
$\SSpof$, in which the weak equivalences are levelwise weak
equivalences of simplicial sets, the fibrations are the levelwise
fibrations of simplicial sets, and the cofibrations are the maps
with the left lifting property with respect to the maps which are
fibrations and weak equivalences.
\end{prop}

Now, we would like to find a map $\varphi_\mathcal O$ with which
to localize $\SSpof$ so that the fibrant objects are Segal
categories.  We first consider the map $\varphi$ used by Rezk
\cite[\S 4]{rezk} to localize simplicial spaces to obtain more
general Segal spaces, then modify it so that the objects are in
$\SSpof$. (Rezk's Segal spaces satisfy the same product condition
as Segal categories. While they are not necessarily discrete in
degree zero, they have the additional requirement that they be
Reedy fibrant.)

Using the maps $\alpha^k$ given above, Rezk defines the object
\[ G(n)^t= \bigcup_{k=0}^{n-1} \alpha^k \Delta [1]^t \] and the inclusion map $\varphi^n:
G(n)^t \rightarrow \Delta [n]^t$.  His localization is with
respect to the coproduct of inclusion maps
\[ \varphi = \coprod_{n \geq 1} (G(n)^t \rightarrow \Delta [n]^t). \]

However, in our case, the objects $G(n)^t$ and $\Delta [n]^t$ are
not necessarily in the category $\SSpo$.  To modify them, we can
replace each $\Delta [n]^t$ with the objects $\Delta [n]^t_\xu$,
where $\xu =(x_0, \ldots ,x_n)$ specifies the 0-simplices. Then,
we define
\[ G(n)^t_\xu = \bigcup_{k=0}^{n-1} \alpha^k \Delta
[1]^t_{x_k, x_{k+1}}. \] Now, we need to take coproducts not only
over all values of $n$, but also over all $n$-tuples of vertices.

So, we define for each $n \geq 0$ the map
\[ \varphi^n = \coprod_{\xu \in \mathcal O^{n+1}}(G(n)^t_\xu
\rightarrow \Delta [n]^t_\xu). \]  Then the map $\varphi$ looks
like
\[ \varphi = \coprod_{n \geq 1}(\varphi^n: \coprod_{\xu \in
\mathcal O^{n+1}} (G(n)^t_\xu \rightarrow \Delta [n]^t_\xu)).
\] When the set $\mathcal O$ is not clear from the context, we
will write $\varphi_\mathcal O$ to specify that we are in $\SSpo$.


For any simplicial space $X$ there is a map
\[ \varphi_n = \Map^h(\varphi^n, X): \Map^h(\coprod_{\xu} \Delta
[n]^t_\xu,X) \rightarrow \Map^h(\coprod_{\xu} G(n)^t_\xu ,X). \]
More simply written, this map is
\[ \varphi_n:X_k \rightarrow \underbrace{X_1 \times_{X_0} \cdots
\times_{X_0} X_1}_n \] and is precisely the Segal map given in the
definition of a Segal category.  Thus, the map $\varphi_\mathcal
O$ is the correct one to use for localizing the model category
$\SSpof$.

\begin{prop} \cite[3.8]{simpmon} \label{sspof}
Localizing the model category structure on $\SSpof$ with respect
to the map $\varphi_\mathcal O$ results in a model category
structure $\LSSpof$ on simplicial spaces with a fixed set
$\mathcal O$ in degree zero in which the weak equivalences are the
$\varphi_\mathcal O$-local equivalences, the cofibrations are
those of $\SSpof$, and the fibrations are the maps with the right
lifting property with respect to the cofibrations which are
$\varphi_\mathcal O$-local equivalences.
\end{prop}

We will find it convenient to make our calculations in a another
model structure with the same weak equivalences but in which all
objects are cofibrant.

\begin{theorem} \cite[3.9]{simpmon} \label{sspoc}
There is a model category structure $\SSpoc$ on the category of
Segal precategories with a fixed set $\mathcal O$ in degree zero
in which the weak equivalences and cofibrations are levelwise, and
in which the fibrations are the maps with the right lifting
property with respect to the acyclic cofibrations. This model
structure can then be localized with respect to the map
$\varphi_\mathcal O$ to obtain a model structure which we denote
$\LSSpoc$.
\end{theorem}


These two model category structures are in fact equivalent to one
another.

\begin{prop} \cite[3.10]{simpmon} \label{equiv}
The adjoint pair given by the identity functor induces a Quillen
equivalence of model categories
\[ \xymatrix@1{\LSSpof \ar@<.5ex>[r] & \LSSpoc. \ar@<.5ex>[l]} \]
\end{prop}

Let $\ast$ denote the set with a single element.  The main result
of \cite{simpmon} is that we have the following chain of Quillen
equivalences, in which the top-most maps are the left adjoints,
which can be composed to form a single Quillen equivalence:
\[ \LSSp_{*,f} \rightleftarrows \LSSetst_* \rightleftarrows \Algtm. \]

\section{A Model for Simplicial Groups}

We would like to extend the results of \cite{simpmon} to the case
of simplicial groups rather than simplicial monoids.  In this
section, we describe the modifications that need to be made in
order to encode the necessary inverses.

In the case of monoids, we consider functors $\Deltaop \rightarrow
\SSets$, where the category $\Deltaop$ has as objects finite
ordered sets $[n]=(0 \rightarrow 1 \rightarrow \cdots \rightarrow
n)$ for each $n \geq 0$ and as morphisms the opposites of the
order-preserving maps between them.  Notice that each $[n]$ can be
regarded as a category with $n+1$ objects and a single morphism $i
\rightarrow j$ whenever $i \leq j$.

Here, we consider instead a category, which we denote $\IDeltaop$,
whose objects are given by small groupoids $I[n]=(0
\rightleftarrows 1 \leftrightarrows \cdots \rightleftarrows n)$
for $n \geq 0$. In other words, each $I[n]$ is a category with
$n+1$ objects and a single isomorphism between any two objects.
The morphisms of $\IDeltaop$are generated by two sets of maps:
first, the opposite of the order-preserving maps, as we have in
$\Deltaop$, and also by a ``flip" morphism on each $I[n]$ which
sends each $i$ to $n-i$. Note that in this case
``order-preserving" should be taken to mean in the standard
numerical ordering of the objects of each $I[n]$, even though
$I[n]$ cannot be considered to be ``ordered" by its morphisms in
the same sense that $[n]$ is.  Alternatively, thinking of
$\IDeltaop$ as a subcategory of the category of the category of
all small groupoids, one should be aware that it is certainly not
a full subcategory.

To understand these maps better, we consider, for example, maps
\[I[2] = (0 \leftrightarrows 1 \leftrightarrows 2) \rightarrow (0 \leftrightarrows 1) =
I[1].\]  Because of the flips, we no longer have only
order-preserving maps, but we do have a preservation of
``betweenness."  Thus, if $0 \mapsto 1$ and $2 \mapsto 1$, then it
follows that $1 \mapsto 1$ also.  In general, if we have a map
$I[n] \rightarrow I[m]$ with $0 \leq i < j < k \leq n$ such that
$i \mapsto \ell$ and $k \mapsto \ell$, for some $0 \leq \ell \leq
m$, it follows that $j \mapsto \ell$ also.

In the case of ${\bf \Delta}$, the simplicial set $\Delta[n]$ is
given by the representable functor $\Hom_{\bf \Delta}(-, [n])$.
Similarly, we can define a simplicial set $I\Delta[n]$ which is
given by the representable functor $\Hom_{\bf I\Delta}(-, I[n])$.
These ``invertible $n$-simplices" are the standard building blocks
of the spaces we consider here.  In particular, every simplex
should be regarded as having a corresponding ``inverse" simplex,
even though this terminology does not make sense in the usual way,
since there is no notion of composition in a simplicial set.  As
with simplicial sets, we can consider the boundary of $I\Delta
[n]$, denoted $I\dot \Delta [n]$, which consists the simplices of
$I \Delta [n]$ of degree less than $n$.

Thus, we can define an \emph{invertible simplicial set} to be a
functor $\IDeltaop \rightarrow \Sets$ and, more generally, an
\emph{invertible simplicial object} in a category $\mathcal C$ to
be a functor $\IDeltaop \rightarrow \mathcal C$.  We denote the
category of invertible simplicial sets by $I\SSets$.  We further
consider the case of invertible simplicial spaces, or functors
$\IDeltaop \rightarrow \SSets$. Since there is a forgetful functor
$U:I\SSets \rightarrow \SSets$ (respectively,
$U:\SSets^{I\Deltaop} \rightarrow \SSets^{\Deltaop}$), we define a
map $f$ of invertible simplicial sets (respectively, spaces) to be
a weak equivalence if $U(f)$ is a weak equivalence of simplicial
sets (respectively, spaces).

In particular, we define a \emph{Segal pregroupoid} to be an
invertible simplicial space $X$ such that the simplicial set $X_0$
is discrete.  A Segal pregroupoid is \emph{reduced} if $X_0 =
\Delta[0]$. To define a Segal groupoid, we need an analogue of the
map $\varphi$ which we used to define Segal categories. However,
just as we defined the maps $\alpha_k:[n]\rightarrow [1]$ in
$\Deltaop$, we can define maps $\beta_k:I[n] \rightarrow I[1]$ in
$\IDeltaop$.  Thus, for any invertible simplicial space $X$ and $n
\geq 2$, we can define the map
\[ \xi_n: X_n \rightarrow \underbrace{X_1 \times_{X_0} \cdots \times_{X_0} X_1}_n \]
induced by the maps
\[ X(\beta_k):X_n \rightarrow X_1. \]  Thus, a \emph{Segal
groupoid} is a Segal pregroupoid $X$ such that for each $n \geq 2$
the map $\xi_n$ is a weak equivalence of invertible simplicial
sets.

As we used the map $\varphi_\mathcal O$ to localize the model
structures $\SSpof$ and $\SSpoc$, we can define an analogous map
$\xi_\mathcal O$ in this situation.  To do so, we first define,
for any $n \geq 2$, the simplicial space $IG(n)^t_\mathcal O$ and
$\xu \in \mathcal O^{n+1}$ given by
\[ IG(n)^t_\mathcal O = \bigcup_{k=0}^{n-1} \beta^k I\Delta[1]^t_{x_k,
x_{k+1}}, \] from which we get an inclusion
\[ \xi^k: \coprod_{\xu \in \mathcal O^{n+1}} (IG(n)^t_\xu \hookrightarrow I
\Delta[n]^t_\xu). \]  Then we have the map
\[ \xi_\mathcal O = \coprod_{n \geq 1}(\xi^n: \coprod_{\xu \in
\mathcal O^{n+1}} (G(n)^t_\xu \rightarrow \Delta [n]^t_\xu)).
\]

\begin{prop}
There is a model category structure $\ISSpof$ on the category of
Segal pregroupoids with a fixed set $\mathcal O$ in degree zero in
which the weak equivalences and fibrations are given levelwise.
Similarly, there is a model category structure $\ISSpoc$ on the
same underlying category in which the weak equivalences and
cofibrations are given levelwise.  Furthermore, we can localize
each of these model category structures with the map $\xi_\mathcal
O$ to obtain model structures $\ILSSpof$ and $\ILSSpoc$ whose
fibrant objects are Segal groupoids.
\end{prop}

\begin{proof}
We begin with the model structure $\ISSpof$, whose existence can
be proved analogously to that of $\SSpof$ \cite[3.7]{simpmon}.  We
first note that, in an abuse of notation, we denote by $K$ the
constant invertible simplicial space given by a simplicial set
$K$.  Furthermore, if $\xu=(x_0, \ldots ,x_n) \in \mathcal
O^{n+1}$, we denote by $I\Delta[n]^t_\xu$ the invertible
simplicial space with the constant simplicial set given by $\xu$
in degree zero.  Specifying the vertices of an $n$-simplex in this
way is necessary since we require all the morphisms to be the
identity on the objects.

Now, we are able to define the objects that we need in order to
define sets of generating cofibrations and generating acyclic
cofibrations.  First, define $(IP_{m,n})_\xu$ to be the pushout of
the diagram
\[ \xymatrix{\dot \Delta[m] \times (I\Delta[n]^t_\xu)_0 \ar[d]
\ar[r] & \dot \Delta[m] \times I \Delta[n]^t_\xu \ar[d] \\
(I\Delta [n]^t_\xu)_0 \ar[r] & (IP_{m,n})_\xu.} \]  Similarly,
define $(IQ_{m,n})_\xu$ to be the pushout
\[ \xymatrix{\Delta[m] \times (I\Delta[n]^t_\xu)_0 \ar[d]
\ar[r] & \Delta[m] \times I \Delta[n]^t_\xu \ar[d] \\
(I\Delta [n]^t_\xu)_0 \ar[r] & (IQ_{m,n})_\xu} \] and
$(IR_{m,n,k})_\xu$ to be the pushout
\[ \xymatrix{V[m,k] \times (I\Delta[n]^t_\xu)_0 \ar[d]
\ar[r] & V[m,k] \times I \Delta[n]^t_\xu \ar[d] \\
(I\Delta [n]^t_\xu)_0 \ar[r] & (IR_{m,n,k})_\xu.} \] Now, we
define the sets
\[ I_f=\{(IP_{m,n})_\xu \rightarrow (IQ_{m,n})_\xu \mid m,n \geq 0, \xu \in \mathcal O^{n+1}\} \] and
\[ J_f= \{(IR_{m,n,k})_\xu \rightarrow (IQ_{m,n})_\xu \mid m,n \geq 0, k \geq 1, \xu \in \mathcal O^{n+1}\}. \]  Then,
applying Theorem \ref{CofGen} to these sets, we can obtain a
cofibrantly generated model category structure $\ISSpof$.

Similarly, by defining appropriate sets of maps, we can obtain
generating cofibrations and generating acyclic cofibrations for
the model structure $\ISSpoc$.  To define these maps, note that
the inclusion functor from the category of Segal pregroupoids to
the category of invertible simplicial spaces has a left adjoint
which can be called a reduction functor, since it ``reduces" the
space in degree zero to a discrete space.  Given an invertible
simplicial space $X$, we denote its reduction by $(X)_r$.

Now, we define sets
\[ I_c=\{(\dot \Delta [m] \times I\Delta[n]^t_\xu \cup \Delta [m]
\times I \dot \Delta [n]^t_\xu)_r \rightarrow (\Delta[m] \times I
\Delta[n]^t_\xu)_r \} \] and
\[ J_c=\{(V[m,k] \times I\Delta[n]^t_\xu \cup \Delta [m] \times I
\dot \Delta[n]^t_\xu)_r \rightarrow (\Delta [m] \times I
\Delta[n]^t_\xu)_r\}.
\] Again, applying Theorem \ref{CofGen}, we obtain the model
structure $\ISSpoc$ just as in \cite[3.9]{simpmon}.

Finally, applying Theorem \ref{localize} with the map
$\xi_\mathcal O$ to each of these model category structures, we
can obtain model category structures $\ILSSpof$ and $\ILSSpoc$.
\end{proof}

As in the Segal category case, Proposition \ref{equiv}, we have
the following result.

\begin{prop}
The adjoint pair given by the identity functor induces a Quillen
equivalence of model categories
\[ \xymatrix@1{\ILSSpof \ar@<.5ex>[r] & \ILSSpoc. \ar@<.5ex>[l]} \]
\end{prop}

We now turn to the theory of groups $\Tg$.  As with the theory of
monoids, we have model structures $\Algtg$ and $\LSSetstg_*$ with
Quillen equivalences
\[ \Algtg \leftrightarrows  \LSSetstg_*,\]
as given in the previous section, where as before $\ast$ denotes
the set with one element.  Thus, we need only show that there is a
Quillen equivalence
\[ \LSSetstg_* \leftrightarrows \ILSSp_{*,f} \] to prove the following
theorem.

\begin{theorem} \label{main}
The model category structure $\Algtg$ is Quillen equivalent to the
model category structure $\ILSSp_{*,f}$.
\end{theorem}

As in \cite{simpmon}, we prove this theorem using several lemmas.
Note that in the model structure $\ILSSp_{*,c}$, we denote by
$L_1$ the localization, or fibrant replacement functor. (Again, it
is convenient to make our calculations in this category rather
than in $\ILSSp_{*,f}$ because here every object is cofibrant.)
Analogously, we denote by $L_2$ the localization functor in
$\LSSetstg_*$.  We will make use of the following general result
for localizations.

\begin{lemma} \cite[4.1]{simpmon} \label{Lhocolim}
Let $L$ be a localization functor on a model category $\mathcal
M$.  Given a small diagram of objects $X_\alpha$ of $\mathcal M$,
\[ L(\hocolim X_\alpha) \simeq L \hocolim (L(X_\alpha)). \]
\end{lemma}

The first step in the proof of the theorem is to show what the
localization functor $L_1$ does to the invertible $n$-simplex
$I\Delta[n]^t$.  By $I\nerve(-)^t$, we denote the invertible nerve
functor $\Hom(I[n],-)$.

\begin{prop} \label{nerve}
Let $F_n$ denote the free group on $n$ generators.  Then in
$\ILSSp_{*, c}$, $L_1 I\Delta [n]^t_*$ is weakly equivalent to
$I\nerve (F_n)^t$ for each $n \geq 0$.
\end{prop}

\begin{proof}
The proof is very similar to the one for reduced Segal categories
\cite[4.2]{simpmon}.  As in that case, note that when $n=0$,
$I\Delta [0]^t_*$ is isomorphic to $I\nerve (F_0)^t$, which is
already a Segal groupoid.

Now we consider the case where $n=1$.  We want to show that the
map $I\Delta [1]^t_* \rightarrow I\nerve(F_1)^t$ obtained by
localizing with respect to the map $\xi$ is a weak equivalence in
$\ILSSp_{*,c}$.

We first define a filtration $\Psi_k$ of $I\nerve(F_1)^t$ whose
set of $j$-simplices looks like
\[ \Psi_k (I\nerve(F_1)^t)_j=\left\{(x^{n_1}|\cdots |x^{n_j}) \mid
\sum_{\ell=1}^j |n_\ell| \leq k \right\} \] where $x$ and its
``inverse" $x^{-1}$ denote the two nondegenerate 1-simplices of $I
\Delta[1]^t_*= \Psi_1$. Thus we have
\[ I\Delta [1]^t_* = \Psi_1 \subseteq \Psi_2 \subseteq
\cdots \subseteq \Psi_k \subseteq \cdots \subseteq \colim_k\Psi_k
.\] We can obtain $\Psi_2$ from $\Psi_1$ by
taking a pushout
\[ \xymatrix{\coprod IG(2)^t_* \ar[r] \ar[d] &
\Psi_1 \ar[d] \\ \coprod I\Delta [2]^t_* \ar[r] & \Psi_2} \] where
the coproducts on the left-hand side are taken over all maps
$IG(2)^t_* \rightarrow \Psi_1$ and $IG(2)^t_*$ is as defined at
the beginning of the section. This process serves to add in, for
example, a ``composite" 1-simplex $x^2$ and the 2-simplex whose
boundary consists of 1-simplices $x$, $x$, and $x^2$, as well as
the inverses to these new simplices. Notice that since we are
working in $\ILSSp_{\ast,c}$, the left-hand vertical map is an
acyclic cofibration, and therefore $\Psi_1 \rightarrow \Psi_2$ is
an acyclic cofibration also \cite[3.14]{ds}.

Similarly, to obtain $\Psi_3$ we will add extra 1-simplices, such
as $x^{-3}$ in order to add a 3-simplex $(x^{-1} | x^{-1} |
x^{-1})$. However, when taking the pushout, we do not want to
start with $IG(3)_*$, since we have already added some of the
1-simplices of this 3-simplex when we localized to obtain
$\Psi_2$.  So, we define $(I\Delta [3]^t_*)_{\Psi_2}$ to be the
part of $I\Delta [3]^t_*$ contained in $\Psi_2$. Then we have a
pushout diagram
\[ \xymatrix{ \coprod (I\Delta [3]^t_*)_{\Psi_2} \ar[r] \ar[d] & \Psi_2 \ar[d] \\
\coprod I\Delta [3]^t_* \ar[r] & \Psi_3} \] where the coproducts
are taken over all maps $(I\Delta[3]^t_*)_{\Psi_2} \rightarrow
{\Psi_2}$. The map $(I\Delta [3]^t_*)_{\Psi_2} \rightarrow I\Delta
[3]^t_*$ is a weak equivalence in $\ILSSp_{\ast,c}$ as follows.
Consider the maps
\[ \xymatrix@1{IG(3)^t_* \ar[r]^-\alpha
& (I\Delta [3]^t_*)_{\Psi_2} \ar[r]^-\beta & I\Delta [3]^t_*}. \]
For any local $X$, the functor $\Map(-,X)$ applied to any of the
three above spaces yields $X_1 \times X_1 \times X_1 \simeq X_3$.
The map $\alpha$ is a weak equivalence since it is just a patching
together of two localizations coming from the map $IG(2)^t_*
\rightarrow I\Delta [2]^t_*$, which is a weak equivalence since it
is one of the maps with respect to which we are localizing. The
composite map $\beta \alpha$ is also a weak equivalence for the
same reason. Thus, $\beta$ is also a weak equivalence by the
two-out-of-three property for weak equivalences. Again, since
$(I\Delta [3]^t_\ast)_{\Psi_2} \rightarrow I\Delta [3]^t_\ast$ is
an acyclic cofibration in $\ILSSp_{\ast, c}$, the map $\Psi_2
\rightarrow \Psi_3$ is an acyclic cofibration also.

For $k>1$, we define $(I\Delta [k+1]^t_*)_{\Psi_k}$ to be the
piece of $I\Delta [k+1]^t_*$ already obtained from previous steps
of the filtration.  Note that it always consists of two copies of
$I\Delta [k]^t_*$ attached along a copy of $I\Delta [k-1]^t_*$, so
the same argument as for $k=2$ shows that the map $(I\Delta
[k+1]^t_*)_{\Psi_k} \rightarrow I\Delta [i+1]^t_*$ is a weak
equivalence.  Hence, for each $k$ we obtain $\Psi_{k+1}$ via the
pushout diagram
\[ \xymatrix{\coprod (\Delta [k+1]^t_*)_{\Psi_k} \ar[r] \ar[d] & \Psi_k \ar[d] \\
\coprod \Delta [k+1]^t_* \ar[r] & \Psi_{k+1}} \] with the
coproducts given as before.

Now that we have defined each stage of our filtration, using the
bar construction notation shows how to map this local object to
$I\text{nerve}(F_1)^t$.  For example,
\[ (x \mid x^{-1} \mid x^2) \mapsto (x, x^{-1}, x^2) \in F_1 \times F_1 \times
F_1.\]

Using Lemma \ref{Lhocolim} we have that
\[ \begin{aligned}
\text{nerve}(F_1)^t & \simeq F_1(I\nerve(F_1)^t) \\
& \simeq L_1(\hocolim (\Psi_k)) \\
& \simeq L_1(\hocolim L_1(\Psi_k)) \\
& \simeq L_1(\hocolim L_1(\Psi_1)) \\
& \simeq L_1L_1(\Psi_1)  \\
& \simeq L_1(\Psi_1) \\
& \simeq L_1(I\Delta [1]^t_*).
\end{aligned} \]

Now, for $n=2$ (i.e., starting with $I\Delta [2]^t_\ast$), we have
six 1-simplices, which we call $x$, $y$, $xy$, $x^{-1}$, $y^{-1}$,
and $(xy)^{-1}$, and two nondegenerate 2-simplices $(x \mid y)$
and $(y^{-1} \mid x^{-1})$. Because we now have two variables, we
need to define the filtration slightly differently as
$\Psi_i=\{[w_1, \ldots ,w_k] \mid \ell(w_1 \ldots w_k) \leq i \}$
where the $w_j$'s are (unreduced) words in $x$ and $y$ and $\ell$
denotes the length of a given word. Note that by beginning with
$\Psi_1$ we start with fewer simplices than those of the 2-simplex
we are considering, but by passing to $\Psi_2$ we obtain $xy$, $(x
\mid y)$ and their respective ``inverses" as well as additional
nondegenerate simplices. (In fact, we are actually starting the
filtration with $\Psi_1=IG(2)^t_*$.) The localizations proceed as
in the case where $n=1$, enabling us to map to $I\nerve (F_2)^t$.

For $n \geq 3$, the same argument works as for $n=2$, with the
filtrations being defined by the lengths of words in $n$ letters.
The resulting object is a reduced Segal category weakly equivalent
to $I\Delta [n]^t_*$. Hence, we have that for any $n$, $L_1
I\Delta [n]^t_*$ is weakly equivalent to $I\nerve (F_n)^t$.
\end{proof}

We now define a functor $J: \IDeltaop \rightarrow \Tg$ induced by
the invertible nerve construction on a group $G$.  For an object
$I[n]$ of ${\bf I\Delta}$, define $J^{op}(I[n]) = F_n$ where $F_n$
denotes the free group on $n$ generators, say $x_1, \ldots ,x_n$.
In particular, $J^{op}(I[0]) = F_0$, the trivial group.  Since $G$
has inverses, there is no difficulty in using the category
$\IDeltaop$ rather than $\Deltaop$.

Taking the (invertible) nerve of a simplicial group $G$ results in
an invertible simplicial space which at level $k$ looks like
\[ I\nerve (G)_k = G^k = \Hom_{\mathcal Grp}(F_k,G). \]
Thus the invertible simplicial diagram $I\nerve (-)$ of
representable functors $\Hom (F_k,-)$ gives rise to an invertible
cosimplicial diagram (i.e., a diagram given by ${\bf I\Delta}$) of
representing objects $T_k$.


To obtain an invertible simplicial diagram of free groups, we
simply reverse the direction of the arrows to obtain a functor $J:
\IDeltaop \rightarrow \Tg$. This map induces a map
$J^*:\SSets^{\Tg} \rightarrow \SSets^{\IDeltaop}$ which can be
restricted to a map $J^*:\SSetstg_* \rightarrow \mathcal
{ISS}p_*$.  Notice that $J^*$ is the identity on objects but
restricts from the morphisms of $\Tg$ to those of $\IDeltaop$.

To obtain a left adjoint to $J^*$, we use a left Kan extension. We
state the following definitions in a general context. Let
$p:\mathcal C \rightarrow \mathcal D$ and $G: \mathcal C
\rightarrow \SSets$ be functors.

\begin{definition}
If $d$ is an object of $\mathcal D$, then the \emph{over category}
or \emph{category of objects over} $d$, denoted $(p \downarrow
d)$, is the category whose objects are pairs $(c,f)$ where $c$ is
an object of $\mathcal C$ and $f:p(c) \rightarrow d$ is a morphism
in $\mathcal D$.  If $c'$ is another object of $\mathcal C$, a
morphism in the over category is given by a map $c \rightarrow c'$
which makes the resulting triangular diagram commute.
\end{definition}

\begin{definition} \cite[11.8.1]{hirsch}
Let $p$, $c$, and $G$ be defined as above, and let $f:p(c)
\rightarrow d$ be an object in $(p \downarrow d)$. The \emph{left
Kan extension} over $p$ is a functor $p_*G: \mathcal D \rightarrow
\SSets$ defined by
\[ (p_*G)(d) = \colim_{(p \downarrow d)}((c, f)
\mapsto G(c)). \]
\end{definition}

Note that, since we are making calculations in the model structure
$\ILSSpoc$ in which all objects are cofibrant, taking the left Kan
extension is homotopy invariant \cite[3.7]{dk3}.

\begin{prop} \cite[11.9.3]{hirsch}
Let $\mathcal C \rightarrow \mathcal D$ be a functor.  The functor
$\SSets^\mathcal C \rightarrow \SSets^\mathcal D$, given by
sending $G$ to the left Kan extension $p_*G$, is left adjoint to
the functor $\SSets^\mathcal D \rightarrow \SSets^\mathcal C$
given by composition with $p$.
\end{prop}

Thus, define $J_*:\mathcal {ISS}p_* \rightarrow \SSets^{\Tg}_*$ to
be the left Kan extension over $J$ which is left adjoint to $J^*$.
Note that even if $G$ is a reduced Segal groupoid, $J_*(G)$ is not
necessarily local in $\LSSetstg_*$. To obtain a $\Tg$-algebra, we
must apply the localization functor $L_2$.  So, we want to know
what we get when we apply $J_*$ followed by $L_2$ to a reduced
Segal groupoid.

Define $IM[k]$ to be the functor $\Tg \rightarrow \SSets$ given by
$F_n \mapsto \Hom_{\Tg}(F_k,F_n)=F_k^n$. Let $H$ be the reduced
Segal groupoid $I\nerve (F_k)^t$.

\begin{lemma}
In $\LSSetstg_*$, $L_2J_*(H)$ is weakly equivalent to $IM[k]$.
\end{lemma}

\begin{proof}
The proof is similar to that of \cite[4.3]{simpmon}. It suffices
to show that for any local object $X$ in $\LSSetstg_*$,
\[\Map^h_{\LSSetstg_*}(L_2J_*H,X) \simeq X(F_k) \] since $\Map^h_{\LSSetstg}(IM[k],X)$ is precisely $X(F_k)$.
This fact can be shown in the following argument:
\[ \begin{aligned}
\Map^h_{\LSSetstg_*}(L_2J_*H,X)& \simeq \Map^h_{\LSSetstg_*}(J_*H,X) \\
& \simeq \Map^h_{\mathcal {ILSS}p_{*,c}}(H,J^*X) \\
& \simeq \Map^h_{\mathcal {ILSS}p_{*,c}}(L_1 I\Delta [k]^t_*, J^*X) \\
& \simeq \Map^h_{\mathcal {ILSS}p_{*,c}}(I\Delta [k]^t_*, J^*X) \\
& \simeq J^*X[k] \\
& \simeq X(F_k).
\end{aligned} \]
\end{proof}

\begin{prop} \label{JLJ}
For any object $X$ in $\mathcal{ISS}p_{*,c}$, we have that $L_1X$
is weakly equivalent to $J^*L_2J_*X$.
\end{prop}

\begin{proof}
This proof is analogous to that of \cite[4.7]{simpmon}. First note
that
\[ X \simeq \hocolim_{\IDeltaop}(I[n] \rightarrow
\coprod_i I\Delta [n_i]^t) \] where the values of $i$ depend on
$n$. We begin by looking at $L_1 X$.  Using Lemma \ref{Lhocolim},
we have the following:
\[ \begin{aligned}
L_1X & \simeq L_1 \hocolim_{\Deltaop} (I[n] \mapsto \coprod
I\Delta [n_i]^t_*)
\\
& \simeq L_1 \hocolim_{\IDeltaop} L_1(I[n] \mapsto \coprod I\Delta [n_i]^t_*) \\
& \simeq L_1 \hocolim_{\IDeltaop} (I[n] \mapsto
\text{nerve}(F_{\sum n_i})^t)
\end{aligned} \]
However, $\hocolim_{\Deltaop} (I[n] \mapsto \nerve(F_{\sum
n_i})^t)$ is already local, a fact which follows from the fact
that the homotopy colimit can be taken at each level, yielding a
Segal pregroupoid in $\mathcal{ILSS}p_{*,c}$ which is still a
Segal groupoid.

Working from the other side of the desired equation, we obtain,
using the fact that left adjoints commute with homotopy colimits:
\[ \begin{aligned}
J^*L_2J_*X & \simeq J^*L_2J_* \hocolim_{\IDeltaop}(I[n] \mapsto
\coprod I\Delta [n_i]^t_*) \\
& \simeq J^*L_2 \hocolim_{\IDeltaop} J_*(I[n] \mapsto \coprod I\Delta [n_i]^t_*) \\
& \simeq J^*L_2 \hocolim_{\IDeltaop} L_2J_*(I[n] \mapsto \coprod I\Delta [n_i]^t_*) \\
& \simeq J^*L_2 \hocolim_{\IDeltaop} (I[n] \mapsto IM[\sum n_i]) \\
\end{aligned} \]

At each level, we have the same spaces as in the
$\IDeltaop$-diagram, but with maps given by $\Tg$ rather than
$\IDeltaop$.  Thus, applying the restriction map $J^*$ results in
a diagram with the same objects at each level, as we wished to
show.
\end{proof}

Combining the above results, the following proposition can be
proved just as \cite[4.8]{simpmon}.

\begin{prop}
The adjoint pair
\[ \xymatrix@1{J_*: \mathcal {ISS}p_{*,f} \ar@<.5ex>[r] & \SSets^{\Tg}_* :J^* \ar@<.5ex>[l]} \]
is a Quillen pair.
\end{prop}

\begin{proof}
We first show that the adjoint pair
\[ \xymatrix@1{J_*: \mathcal {SS}p_{*,f} \ar@<.5ex>[r] & \SSets^{\Tm}_* :J^* \ar@<.5ex>[l]} \]
is a Quillen pair.  In both model categories, the fibrations and
weak equivalences are defined levelwise.  Since the right adjoint
$J^*$ preserves the spaces at each level, it must preserve both
fibrations and acyclic fibrations.

An application of Theorem \ref{LocPair} shows that we still have a
Quillen pair after the respective localizations.
\end{proof}

\begin{proof}[Proof of Theorem \ref{main}]
First, we need to know that the right adjoint $J^*$ reflects weak
equivalences between fibrant objects.  In each of the two
localized model categories $\ILSSp_{\ast, f}$ and $\LSSetstg_*$,
an object is fibrant if and only if it is local and fibrant in the
unlocalized model category. Therefore, in each case a weak
equivalence between fibrant objects is a levelwise weak
equivalence.  Since $J^*$ does not change the spaces at each
level, it must reflect weak equivalences between fibrant objects.

Finally, by Proposition \ref{JLJ}, $L_1X \simeq J^*L_2J_*X$ for
any functor $X: \Deltaop \rightarrow \SSets$, and in particular
for any cofibrant $X$.
\end{proof}

Now composing the Quillen equivalences
\[ \ILSSp_{*,f} \rightleftarrows \LSSetstg_* \rightleftarrows
\Algtg \] (where the left adjoint functors are the topmost maps)
results in a single Quillen equivalence
\[ \ILSSp_{*,f} \rightleftarrows \Algtg. \]

\section{Segal Groupoids}

In this section, we generalize the result on simplicial groups and
reduced Segal pregroupoids to one on simplicial groupoids and
Segal groupoids.  While Segal groupoids (and their $n$-categorical
analogues) have been discussed in works by Hirshowitz and Simpson
\cite{hs}, \cite{simpson}, we give an alternate but equivalent
description of them.

\begin{definition}
A \emph{simplicial category} is a category enriched over
simplicial sets, or a category in which there is a simplicial set
of morphisms between any two objects.  A \emph{simplicial
groupoid} is a simplicial category enriched over invertible
simplicial sets.
\end{definition}


For convenience, we describe a modification of the categories
$\Deltaop$ and $I\Deltaop$.  Let $\mathcal O$ be a set. We define
the category $\Deltaop_{\mathcal O}$ as follows. The objects are
given by $[n]_{x_0, \ldots ,x_n}$ where $n \geq 0$ and $(x_0,
\ldots x_n) \in \mathcal O^{n+1}$. The $[n]$ should be thought of
as in the simplicial category $\Deltaop$; however, recall that
when we work with Segal categories we will require all morphisms
to preserve the objects. Therefore, we need to have a separate
$[n]$ for each possible $(n+1)$-tuple of objects in $\mathcal O$.
The morphisms in $\Deltaop_\mathcal O$ are those of $\Deltaop$ but
depend on the choice of $(x_0, \ldots ,x_n)$. Specifically, the
face maps are
\[ d_i:[n]_{x_0, \ldots ,x_n} \rightarrow [n-1]_{x_0, \ldots
,\widehat x_i, \ldots ,x_n} \] and the degeneracy maps are
\[ s_i:[n]_{x_0, \ldots ,x_n} \rightarrow [n+1]_{x_0, \ldots,
x_{i-1}, x_i, x_i, x_{i+1}, \ldots ,x_n}. \]  Note that if
$\mathcal O$ is the one-object set, then $\Deltaop_{\mathcal O}$
is just $\Deltaop$.

Similarly, we can describe the category $\IDeltaop_\mathcal O$.
Note that in this case the flip map is no longer a self-map
(unless $\mathcal O$ consists of a single element), but instead a
map $[n]_{x_0, \ldots, x_n} \rightarrow [n]_{x_n, \ldots,x_0}$.

Now we can use this notation to describe Segal categories and
Segal groupoids. A Segal category with $\mathcal O$ in degree zero
is a functor $X:\Deltaop_\mathcal O \rightarrow \SSets$ such that
for each $n \geq 2$ and $(n+1)$-tuple $(x_0, \ldots ,x_n)$, the
map
\[ X([n]_{x_0, \ldots, x_n}) \rightarrow X([1]_{x_0, x_1}) \times
_{X[0]_{x_1}} \cdots \times_{X[0]_{x_{n-1}}} X([1]_{x_{n-1},x_n})
\] is a weak equivalence.  Analogously, a Segal groupoid with
$\mathcal O$ in degree zero is a functor $X:I\Deltaop_\mathcal O
\rightarrow \SSets$ satisfying these same conditions.

Recall from section 2 that we have model category structures
$\LSSpof$ and $\LSSpoc$ on the category of Segal precategories
with $\mathcal O$ in degree zero, in each of which the fibrant
objects are Segal categories.  Then, from section 3, we have the
analogous model structures $\ILSSpof$ and $\ILSSpoc$, whose
fibrant objects are Segal groupoids.

We would like to think of the category of simplicial categories,
or the category of simplicial groupoids, with object set $\mathcal
O$ as a diagram category as well. To do so, we need to define the
notion of a multi-sorted algebraic theory. To see more details,
see \cite{multisort}.

\begin{definition}
Given a set $S$, an $S$-\emph{sorted algebraic theory} (or
\emph{multi-sorted theory}) $\mathcal T$ is a small category with
objects $T_{\alphau^n}$ where $\alphau^n = <\alpha_1, \ldots
,\alpha_n>$ for $\alpha_i \in S$ and $n \geq 0$ varying, and such
that each $T_{\alphau^n}$ is equipped with an isomorphism
\[ T_{\alphau^n} \cong \prod_{i=1}^n T_{\alpha_i}. \]
For a particular $\alphau^n$, the entries $\alpha_i$ can repeat,
but they are not ordered.  There exists a terminal object $T_0$
(corresponding to the empty object of $S$).
\end{definition}

In particular, we can talk about the theory of $\mathcal
O$-categories, which we will denote by $\Tocat$, and the theory of
$\mathcal O$-groupoids, which we denote $\Togd$. To define these
theories, first consider the category $\mathcal {OC}at$ whose
objects are the categories with a fixed object set $\mathcal O$
and whose morphisms are the functors which are the identity map on
the objects. The objects of $\Tocat$ are categories which are
freely generated by directed graphs with vertices corresponding to
the elements of the set $\mathcal O$. This theory is sorted by
pairs of elements in $\mathcal O$, corresponding to the morphisms
with source the first element and target the second.  In other
words, this theory is $(\mathcal O \times \mathcal O)$-sorted
\cite[3.5]{multisort}. (In the one-object case, we get the
ordinary theory of monoids, since a monoid is just a category with
one object.) We can then say that a simplicial category with
object set $\mathcal O$ is essentially a strict $\Tocat$-algebra,
where the definitions of strict and homotopy $\mathcal T$-algebras
for multi-sorted theories $\mathcal T$ are defined analogously to
those for ordinary algebraic theories. The objects of $\Togd$ are
representatives of the isomorphism classes of finitely generated
free groupoids with object set $\mathcal O$, and a simplicial
groupoid is essentially a strict $\Togd$-algebra.

Again, we have a model structure $\Algt$ on the category of all
$\mathcal T$-algebras and a model structure $\SSetst$ on the
category of all functors $\mathcal T \rightarrow \SSets$ which can
be localized as before to obtain a model category structure
$\LSSetst$ in which the local objects are homotopy $\mathcal
T$-algebras \cite[4.11]{multisort}.  For $\Tocat$, we can define a
category $\SSetsTocat_\mathcal O$ of functors $\Tocat \rightarrow
\SSets$ which send $T_0$ to $\coprod_\mathcal O \Delta [0]$.
Making modifications as in the case of $\LSSetstm_*$, we can
define a model structure $\LSSetsTocat_\mathcal O$ which is
Quillen equivalent to $\Algtocat$ \cite{simpmon2}. Similarly, we
can define $\LSSetsTogd_\mathcal O$ and show that it is Quillen
equivalent to $\Algtogd$.

In particular, Theorem \ref{rigid} for algebraic theories can be
generalized to the case of multi-sorted theories, and therefore we
have that there is a Quillen equivalence of model categories
between $\Algt$ and $\LSSetst$ for any multi-sorted theory
$\mathcal T$ \cite[5.1]{multisort}. We can again use the version
with the stricter requirement on degree zero to obtain a Quillen
equivalence
\[ \LSSetsTocat_\mathcal O \rightleftarrows \mathcal Alg^{\Tocat}.
\]  Hence, the problem reduces to finding a Quillen equivalence between
$\LSSpof$ and $\mathcal{LSS}ets^{\Tocat}_\mathcal O$.

\begin{theorem} \cite[5.5]{simpmon}
The adjoint pair
\[ \xymatrix@1{J_*: \LSSpof \ar@<.5ex>[r] & \LSSetsTocat_\mathcal O :J^* \ar@<.5ex>[l]} \]
is a Quillen equivalence.
\end{theorem}

We can extend this result just as we did in the previous section
to obtain a result relating Segal groupoids with object set
$\mathcal O$ to homotopy algebras over a multi-sorted theory of
groupoids with object set $\mathcal O$.

The same methods can be used to prove the following.

\begin{theorem}
There exists an adjoint pair
\[ \xymatrix@1{J_*: \ILSSpof \ar@<.5ex>[r] & \LSSetsTogd_\mathcal O :J^* \ar@<.5ex>[l]} \]
which is a Quillen equivalence.
\end{theorem}

Hence, composing Quillen equivalences results in a Quillen
equivalence
\[ \Algtgd \leftrightarrows \ILSSpof. \]

The more general result that there is a Quillen equivalence
between model structures on the category of all Segal
precategories and the category of all small simplicial categories
\cite[8.6]{thesis} can be extended to the groupoid situation, but
we defer its proof to the next paper \cite{groupoid}.

\section{An Alternative Model for Simplicial Groups}

In this section, we summarize a result of Bousfield for modelling
simplicial groups with ${\bf \Delta}$ by changing the projection
maps and therefore the Segal maps \cite{bous}.  This approach is
more convenient for adaptation to the case of abelian groups, as
we show in the next section.


Let us recall the way we used ${\bf \Delta}$ to obtain a model for
simplicial monoids.  Given a reduced Segal precategory $X:\Deltaop
\rightarrow \SSets$, the Segal map $X_n \rightarrow (X_1)^n$ was
induced by $\alpha^k:[1] \rightarrow [n]$ in ${\bf \Delta}$ for
each $0 \leq k \leq n-1$, where $\alpha^k(0)=k$ and
$\alpha^k(1)=k+1$.  Heuristically, we are thinking of the
$(X_1)^n$ as a chain of $k$ morphisms which has a ``composite" (at
least up to homotopy) if the Segal map is to be a weak
equivalence.

To consider a basic example, consider a 2-simplex with vertices
$\{a,b,c\}$ and 1-simplices $\{a \rightarrow b, b\rightarrow c, a
\rightarrow c\}$.  The 2-chain here is $a \rightarrow b
\rightarrow c$, the composite is given by $a \rightarrow c$, and
the 2-simplex is defined by all three.  The first projection map
sends the 2-simplex to $a \rightarrow b$ and the second projection
sends it to $b \rightarrow c$.

The idea behind Bousfield's construction is to define these
projections differently.  In the situation just described, the
first projection remains the same.  The second projection,
however, sends the 2-simplex to $a \rightarrow c$.  Thus, if we
are going to fill in the third 1-simplex, we get the ``inverse" of
the first projection composed with the second projection.

To formalize this construction, we define in ${\bf \Delta}$ the
maps $\gamma^k:[1] \rightarrow [n]$ given by $0 \mapsto 0$ and $1
\mapsto k+1$ for all $0 \leq k < n$.  Again restricting to Segal
precategories, we can consider the Bousfield-Segal map $\psi_n:X_n
\rightarrow (X_1)^n$ induced by these maps. The models for
simplicial groups in this sense will be the reduced Segal
categories for which the Bousfield-Segal maps are weak
equivalences of simplicial sets for all $n \geq 2$.  We call such
simplicial spaces \emph{reduced Bousfield-Segal categories}.

To give a localized model structure, we define for each $k \geq 2$
the simplicial space
\[ H(k)^t_* = \bigcup_{i=1}^{k-1} \gamma^i \Delta[1]^t_* \subseteq
\Delta[k]^t_*. \]  Then, as in the previous situations, define the
map
\[ \psi_*= \coprod_{k \geq 1}(\psi^k: H(k)^t_* \rightarrow \Delta [k]^t_*).
\]

\begin{prop}
Localizing the model category structure $\SSp_{*,f}$ with respect
to the map $\psi_*$ results in a model category structure
$\mathcal L_B \SSp_{*,f}$ whose fibrant objects are reduced
Bousfield-Segal categories.  There is also an analogous model
structure $\mathcal L_B \SSp_{*,c}$.
\end{prop}

We denote by $L_B$ the localization functor in $\mathcal L_B
\SSpoc$. With some minor technical changes, our proof that Segal
categories model simplicial monoids generalizes to show that
Bousfield-Segal categories model simplicial groups. As with the
previous proof, the key point is the analogue of Lemma
\ref{nerve}.

\begin{prop}
Let $F_n$ denote the free group on $n$ generators.  Then in
$\mathcal L_B \SSpoc$, $L_B\Delta[n]^t_*$ is weakly equivalent to
$\nerve(F_n)^t$.
\end{prop}

\begin{proof}
In the case where $n=0$, we have that $\Delta[0]^t_* \cong
\nerve(F_0)^t$.  So, we consider the case where $n=1$.  As before,
we define a filtration $\Psi_1 \subseteq \Psi_2 \subseteq \cdots
\subseteq \Psi_k \subseteq \cdots$  We will use a bar construction
notation as before, but, as we are assuming each slot is given by
the image of a projection map, its meaning has changed as we have
changed our projection maps.  The set of $j$-simplices (for $j>1$)
of $\Psi_k$ is given by
\[ \Psi_k (\text{nerve}(F_1)^t)_j=\{(x^{n_1}|\cdots |x^{n_j}) \} \]
with the following additional conditions on the superscripts
$n_\ell$.  As before, we require
\[ \sum_{\ell=1}^j |n_\ell| \leq k. \]  Furthermore, the construction
imposes the following conditions depending on $k$.  For
$k=1,2$, we have $0 \leq n_\ell \leq 1$.  For $k=3$, we have $-1
\leq n_\ell \leq 1$.  For all $k \geq 4$, we have $-k+1 \leq
n_\ell \leq k-2$. Note that our definition of $\Psi_1$ coincides
with that of $\Delta[1]^t_*$ since each has only one
non-degenerate 1-simplex $x$ and no nondegenerate 2-simplices.
Thus we have
\[ \Delta [1]^t_* = \Psi_1 \subseteq \Psi_2 \subseteq
\cdots \subseteq \Psi_k \subseteq \cdots \colim_m \Psi_m. \]

Setting $k=2$, we can take a pushout
\[ \xymatrix{\coprod H(2)^t_* \ar[r] \ar[d] &
\Psi_1 \ar[d] \\ \coprod \Delta [2]^t_* \ar[r] & \Psi_2} \] where
the coproducts on the left-hand side are over the maps $H(2)^t_*
\rightarrow \Psi_1$.  In doing so, we obtain the 1-simplex
$x^{-1}$ and 2-simplices given by $(1 \mid 1)$, $(1 \mid x)$, $(x
\mid 1)$, and $(x \mid x)$, where the third is the one which
requires $x^{-1}$. (Note that our use of the bar construction
notation has been adapted to our new projections.)

Then for $k \geq 3$, $\Psi_k$ is obtained from $\Psi_{k-1}$ by
taking pushouts along the maps $\coprod
(\Delta[k]^t_*)_{\Psi_{k-1}} \rightarrow \coprod \Delta[k]^t_*$
for each $k \geq 3$, where $(\Delta[k]^t_*)_{\Psi_{k-1}}$ denotes
the piece of $\Delta^t[k]_*$ which we have already obtained.  As
usual, the coproduct is taken over all maps
$(\Delta[k]^t_*)_{\Psi_{k-1}} \rightarrow \Psi_{k-1}$.

Then, we can use the inclusions
\[ H(k)^t_* \rightarrow (\Delta[k]^t_*)_{\Psi_{k-1}} \rightarrow
\Delta[k]^t_* \] to prove the fact that each map $\Psi_{k-1}
\rightarrow \Psi_k$ is a weak equivalence, as in the proof of
Lemma \ref{nerve} or that of \cite[4.2]{simpmon}.
\end{proof}

The rest of the proof follows similarly to the one in section 4.
As in section 5, it can also be generalized to obtain a result for
more general Segal groupoids.  We should note that Bousfield's
result was stated in terms of homotopy categories and that his
proof did not make use of model categories.  Furthermore, he
worked in more generality, finding a model for $n$-fold loop
spaces.

\section{A Model for Simplicial Abelian Groups}

In \cite{segal}, Segal defines $\Gamma$-spaces and shows that
strict $\Gamma$-spaces are equivalent to simplicial abelian
monoids.  Here, we begin with his definition of the category
$\Gamma$.  Its objects are representatives of isomorphism classes
of finite sets, and a morphism $S \rightarrow T$ is given by a map
$\theta:S \rightarrow \mathcal P(T)$ such that $\theta(\alpha)$
and $\theta(\beta)$ are disjoint whenever $\alpha \neq \beta$. (We
denote by $\mathcal P(T)$ the power set, or set of all subsets of
the set $T$.)

We can then define the opposite category $\Gammaop$, which has the
following description of its own.  It is the category with objects
${\bf n}=\{0, 1, \ldots, n\}$ for $n \geq 0$ and morphisms ${\bf
m} \rightarrow {\bf n}$ such that $0 \mapsto 0$.  Here, we find it
convenient to use both descriptions, depending on the situation.

Segal defines a $\Gamma$-space $X$ to be a functor $\Gamma
\rightarrow \SSets$ such that $X_0 \simeq \Delta[0]$ and the Segal
map $\varphi_k:X_k \rightarrow (X_1)^k$ is a weak equivalence of
simplicial sets.  He further mentions that if $X_0=\Delta[0]$ and
if each Segal map $\varphi_k:X_k \rightarrow (X_1)^k$ is an
isomorphism rather than a weak equivalence, i.e., if $X$ is a
\emph{strict} $\Gamma$-\emph{space}, then $X$ is essentially a
simplicial abelian monoid.  (The following proof sketch is due to
Badzioch.)

\begin{prop} \cite{segal}
The category of strict $\Gamma$-spaces is equivalent to the
category of simplicial abelian monoids.
\end{prop}

\begin{proof}[Sketch of proof]
First, recall that the category of simplicial abelian monoids is
equivalent the the category $\Algtam$, where $\Tam$ is the theory
of abelian monoids.

Among the maps of $\Gammaop$, there are projections $p_{n,i}:{\bf
n} \rightarrow {\bf 1}$, where $p_{n,i}(k)=1$ if $k=i$ and 0
otherwise.  Then there exists a projection-preserving functor
$\Gammaop \rightarrow \Tam$, where $\Tam$ denotes the theory of
abelian monoids.  Given a strict $\Gamma$-space $X$, it is
uniquely determined by each $X_n$, the projection maps $X_n
\rightarrow X_1$, and the map $X_2 \rightarrow X_1$ which is the
image of the map ${\bf 2} \rightarrow {\bf 1}$ given by $0 \mapsto
0$ and $1,2 \mapsto 1$. In particular, this map $X_2 \rightarrow
X_1$ induces (by induction) all maps $X_n \rightarrow X_1$ arising
from the morphisms ${\bf n} \rightarrow {\bf 1}$ given by $0
\mapsto 0$ and $i \mapsto 1$ for all $0<i \leq n$.

Then, a strict $\Gamma$-space gives the space $X_1$ the structure
of an abelian monoid with multiplication map given by the
specified map $X_2 \rightarrow X_1$.  In particular, $X_1$ defines
a $\Tam$-algebra $tX$ which is given by the necessary $X_n$ for
each $n \geq 2$ and the projection maps. Then, notice that
restricting the $\Tam$-algebra $tX$ to $\Gammaop$ results in our
original $X$.  In other words, if $F$ is the forgetful functor
from the category of $\Tam$-algebras to the category of strict
$\Gammaop$-spaces, we have that $tF(X)=X$. Thus, the functors $t$
and $X$ are inverse to one another.
\end{proof}

Segal defines a functor ${\bf \Delta} \rightarrow \Gamma$ as
follows.  The object $[n]$ is sent to ${\bf n}$ for each $n \geq
0$, and a map $f \colon [m] \rightarrow [n]$ is sent to the map
$\theta \colon {\bf m} \rightarrow {\bf n}$ given by
$\theta(i)=\{j \in {\bf n} \mid f(i-1)<j \leq f(i)\}$.  In
particular, the maps $\alpha^k \colon [1] \rightarrow [n]$ are
sent to maps $\theta \colon {\bf 1} \rightarrow {\bf k}$ given by
$\theta(1)=\{k+1\}$.

Now, we'd like to know what happens if we define the projections
as Bousfield does.  Given $\gamma^k \colon [1] \rightarrow [n]$ in
${\bf \Delta}$, it is sent to the map $j^k \colon {\bf 1}
\rightarrow {\bf n}$ given by $j^k(1)=\{1, \ldots ,k+1\}$. Thus,
we can use the maps $j_k \colon {\bf 1} \rightarrow {\bf n}$ to
define a modified version of the conditions for a $\Gamma$-space.

Now, a \emph{strict Bousfield} $\Gamma$-\emph{space} is a functor
$X:\Gamma^{op} \rightarrow \SSets$ such that $X({\bf 0})=
\Delta[0]$ and the maps $X({\bf n}) \rightarrow X({\bf 1})^n$
induced by the maps $j_k$ for all $1 \leq k \leq n$ are
isomorphisms for all $n \geq 2$. Similarly, a \emph{(homotopy)
Bousfield} $\Gamma$-\emph{space} has $X({\bf 0})$ contractible and
the above map a weak equivalence of simplicial sets.  The above
proof for strict $\Gamma$-spaces extends to this new situation.

\begin{theorem}
There is an equivalence between the category of strict Bousfield
$\Gamma$-spaces and the category of simplicial abelian groups.
\end{theorem}

However, it is expected, in analogy with the original work of
Segal \cite{segal}, that the (homotopy) Bousfield $\Gamma$-spaces
are not equivalent to simplicial abelian groups up to homotopy.
Since Segal shows that $\Gamma$-spaces are equivalent to infinite
loop spaces when they have homotopy inverses \cite[1.4]{segal}, it
is expected that the Bousfield $\Gamma$-spaces are equivalent to
infinite loop spaces.  Another interesting comparision to be made
is with the work of Schw\"{a}nzl and Vogt \cite{sv}.

\end{document}


\title[Erratum to ``Adding inverses"]{Erratum to ``Adding Inverses to Diagrams Encoding Algebraic Structures" and ``Adding inverses to diagrams II: Invertible homotopy theories are spaces"}

\author{Julia E. Bergner}
\email{bergnerj@member.ams.org}

\address{University of California, Riverside, CA 92521}

\subjclass[2010]{55U35, 18G30, 18E35}

\thanks{The author was partially supported by NSF grant DMS-0805951.}

\maketitle

\section{Statement of previous error}

In previous work, we studied various kinds of functors $X \colon \Deltaop \rightarrow \SSets$ satisfying a Segal condition, so that the maps
\[ X_n \rightarrow \underbrace{X_1 \times_{X_0} \cdots \times_{X_0} X_1}_n \] are weak equivalences of simplicial sets for $n \geq 2$.  When we imposed the additional condition that $X_0=\Delta[0]$, such objects, called \emph{reduced Segal categories} or \emph{Segal monoids} were shown to be equivalent to simplicial monoids \cite{simpmon}.  When instead $X_0 =\amalg \Delta[0]$, some discrete simplicial set, then $X$ is a \emph{Segal category} can more generally be regarded as an up-to-homotopy model for a simplicial category with this same object set \cite{simpmon}.

These results were used in the comparison between the model structure for all Segal categories (not just with a fixed set in degree zero) and the model structure for simplicial categories \cite{thesis}.  Furthermore, the former model structure was shown to be Quillen equivalent to the model structure for complete Segal spaces, or functors $\Deltaop \rightarrow \SSets$ where the discrete level zero condition is replaced with a ``completeness" condition \cite{thesis}.

A natural question was then whether these results could be generalized to an ``invertible" and two methods were proposed in both of the papers \cite{inverses}, \cite{inverses2}.  The first was to replace $\Deltaop$ in the above definitions with a category $\IDeltaop$ in which the objects had an involution map.  However, these results were in fact incorrect, in that this involution does not adequately encode an inverse map. In this note, we clarify that this diagram should encode the structure of a monoid with involution rather  than a group, or category with involution rather than a groupoid, in the case of multiple objects.

The second approach given in these papers, given by using different projection maps as first used by Bousfield \cite{bous}, is still correct.

\begin{thank}
The author would like to thank Philip Hackney for conversations about this work.
\end{thank}

\section{Replacing groups with monoids with involution}

In the case of monoids, we consider functors $\Deltaop \rightarrow
\SSets$, where the category $\Deltaop$ has as objects finite
ordered sets $[n]=(0 \rightarrow 1 \rightarrow \cdots \rightarrow
n)$ for each $n \geq 0$ and as morphisms the opposites of the
order-preserving maps between them.  Notice that each $[n]$ can be
regarded as a category with $n+1$ objects and a single morphism $i
\rightarrow j$ whenever $i \leq j$.

In \cite{inverses} we defined a category $\IDeltaop$,
whose objects are given by small groupoids $I[n]=(0
\rightleftarrows 1 \leftrightarrows \cdots \rightleftarrows n)$
for $n \geq 0$. In other words, each $I[n]$ is a category with
$n+1$ objects and a single isomorphism between any two objects.
The morphisms of $\IDeltaop$ are generated by two sets of maps:
the opposite of the order-preserving maps from
$\Deltaop$, and an involution morphism on each $I[n]$ which
sends each $i$ to $n-i$.

The hope was that functors $\IDeltaop \rightarrow \SSets$ satisfying a Segal condition encoded a group structure.  Unfortunately, inverses are not adequately given, so such functors actually give the structure of a monoid with involution.  Thus, we change the terminology given in \cite{inverses} as follows.

In the case of ${\bf \Delta}$, the simplicial set $\Delta[n]$ is
given by the representable functor $\Hom_{\bf \Delta}(-, [n])$.
Similarly, we can define an object $I\Delta[n]$ which is
given by the representable functor $\Hom_{\bf I\Delta}(-, I[n])$.
These $n$-\emph{simplices with involution} are the standard building blocks
of the spaces we consider here.  In particular, every simplex
should be regarded as having a corresponding ``reverse" simplex.  As
with simplicial sets, we can consider the boundary of $I\Delta
[n]$, denoted $\partial I \Delta [n]$, which consists of the nondegenerate simplices of
$I \Delta [n]$ of degree less than $n$.

Thus, we can define an \emph{simplicial set with involution} to be a
functor $\IDeltaop \rightarrow \Sets$ and, more generally, an
\emph{simplicial object with involution} in a category $\mathcal C$ to
be a functor $\IDeltaop \rightarrow \mathcal C$.  We denote the
category of simplicial sets with involution by $I\SSets$.  We further
consider the case of simplicial spaces with involution, or functors
$\IDeltaop \rightarrow \SSets$. Since there is a forgetful functor
$U:I\SSets \rightarrow \SSets$ (respectively,
$U:\SSets^{I\Deltaop} \rightarrow \SSets^{\Deltaop}$), we define a
map $f$ of simplicial sets (respectively, spaces) with involution to be
a weak equivalence if $U(f)$ is a weak equivalence of simplicial
sets (respectively, spaces).

In particular, we define a \emph{Segal precategory with involution} to be an simplicial space with involution $X$ such that the simplicial set $X_0$
is discrete.  If $X_0 =\Delta[0]$, then we call it a \emph{Segal premonoid with involution}. To define a Segal category with involution, use the maps
\[ \xi_n: X_n \rightarrow \underbrace{X_1 \times_{X_0} \cdots \times_{X_0} X_1}_n \] defined in \cite[\S 4]{inverses}. Thus, a \emph{Segal category with involution} is a Segal precategory with involution $X$ such that for each $n \geq 2$
the map $\xi_n$ is a weak equivalence of simplicial sets.

Obtaining an appropriate model structure requires localization with respect to the following map:
\[ \xi_\mathcal O = \coprod_{n \geq 1}(\xi^n: \coprod_{\xu \in \mathcal O^{n+1}} (IG(n)^t_\xu \rightarrow I\Delta [n]^t_\xu)). \]

The proofs of the following two propositions continue to hold, with the necessary changes in terminology.

\begin{prop} \cite[4.1]{inverses}
There is a model category structure $\LSSets^{\IDeltaop}_{\mathcal O,f}$ on the category of
Segal precategories with involution with a fixed set $\mathcal O$ in degree zero in
which the weak equivalences and fibrations are given levelwise.
Similarly, there is a model category structure $\LSSets^{\IDeltaop}_{\mathcal O,c}$ on the
same underlying category in which the weak equivalences and
cofibrations are given levelwise.  Furthermore, we can localize
each of these model category structures with the map $\xi_\mathcal
O$ to obtain model structures $\LSSets^{\IDeltaop}_{\mathcal O,f}$ and $\LSSets^{\IDeltaop}_{\mathcal O,c}$ whose
fibrant objects are Segal categories with involution.
\end{prop}

\begin{prop} \cite[4.2]{inverses}
The adjoint pair given by the identity functor induces a Quillen
equivalence of model categories
\[ \xymatrix@1{\LSSets^{\IDeltaop}_{\mathcal O,f} \ar@<.5ex>[r] & \LSSets^{\IDeltaop}_{\mathcal O,c}. \ar@<.5ex>[l]} \]
\end{prop}

In \cite{inverses}, we claimed that there was a Quillen equivalence
\[ \LSSetstg_* \leftrightarrows \LSSets^{\IDeltaop}_{*,f}. \]  Unfortunately, we did not adequately establish that we obtained group structures using this category $\IDeltaop$.  We seek to establish that our previous proof instead gave a Quillen equivalence
\[ \LSSetstmi_* \leftrightarrows \LSSets^{\IDeltaop}_{*,f} \] where $\Tmi$ is the \emph{theory of monoids with involution}.

This theory $\Tmi$ has as objects $T_k$ which are given by the free monoid with involution on $k$ generators.  In other words, $T_k$ is free on generators $x_1, \ldots, x_k$ and $\overline{x}_1, \ldots, \overline{x}_k$ with involution $I(x_k)=\overline{x}_k$.  There is a monoid map $\tau \colon T_k \rightarrow T_k$ given by $x_j \mapsto \overline{x}_{k-j+1}$.  This map will correspond to the flip map of each object $I[k]$ in $\IDeltaop$.

Using Badzioch's theorem from \cite{bad}, such a Quillen equivalence will complete the proof of the following theorem.

\begin{theorem} \label{main}
The model category structure $\Algtmi$ is Quillen equivalent to the
model category structure $\LSSets^{\IDeltaop}_{*,f}$.
\end{theorem}

As in \cite{simpmon}, we prove this theorem using several lemmas.
Note that in the model structure $\LSSets^{\IDeltaop}_{*,c}$, we denote by
$L_1$ the localization, or fibrant replacement functor.
Analogously, we denote by $L_2$ the localization functor in
$\LSSetstmi_*$.

The first step in the proof of the theorem is to show what the
localization functor $L_1$ does to the $n$-simplex with involution
$I\Delta[n]^t$.  By $I\nerve(-)^t$, we denote the representable
functor $\Hom(I[n],-)$, viewed as a transposed constant simplicial space.  It is here, Lemma \cite[4.5]{inverses}, that the major error occurred.  The correct statement is as follows.

\begin{prop}  \label{nerve}
Let $F_n$ denote the free monoid with involution on $n$ generators.  Then in
$\LSSets^{\IDeltaop}_{*,c}$, $L_1 I\Delta [n]^t_*$ is weakly equivalent to $I\nerve (F_n)^t$ for each $n \geq 0$.
\end{prop}

In the proof, we defined a filtration of $I\nerve(F_1)^t$ as follows:
\[ \Psi_k (I\nerve(F_1)^t)_j=\left\{(x^{n_1}|\cdots |x^{n_j}) \mid \sum_{\ell=1}^j |n_\ell| \leq k \right\} \]
where $x$ and its ``inverse" $x^{-1}$ denote the two nondegenerate 1-simplices of $I \Delta[1]^t_*= \Psi_1$.  The problem is that we assumed here that $x$ could be canceled with $x^{-1}$, when there is no structure built in to $\IDeltaop$ to make this cancelation possible.  In short, we do not really have the structure of a group, but only of a monoid with involution.  Therefore, we instead define the filtration of $I\nerve(F_1)^t$ by
\[ \Psi_k (I\nerve(F_1)^t)_j=\left\{(w_1|\cdots |w_j) \mid \ell(w_1 \cdots w_j) \leq k \right\} \] where the $w_i$ are words in $x$ and $\overline{x}$, and $\ell$ denotes word length.

With this modification, the previous proof goes through as before.

To obtain the appropriate adjoint pair for the Quillen equivalence, we first define a functor $J \colon \IDeltaop \rightarrow \Tmi$.  On objects, this functor is defined by $I{n} \mapsto T_n$.  On face and degeneracy maps coming from $\Deltaop$, this functor behaves the same as the functor $\Deltaop \rightarrow \Tm$ as defined in \cite{simpmon}.  Specifically, the coface and codegeneracy maps maps in ${\bf I\Delta}$ are given by
\[ d^i(x_k)=
\begin{cases}
x_k & i>k \\
x_kx_{k+1} & i=k \\
x_{k+1} & i<k
\end{cases} \text{ and }
s^i(x_k)=
\begin{cases}
x_k & i \leq k \\
e & i=k-1 \\
x_{k-1} & i<k-1.
\end{cases} \]
The involution map is sent to the map $\tau \colon T_\ell \rightarrow T_\ell$ defined by $x_k \mapsto y_{\ell-k+1}$.

This functor $J$ induces a map
\[ J^\ast \colon \LSSetstmi_{\ast,f} \rightarrow \LSSets^{\IDeltaop}_{\ast,f} \] for which we have a left adjoint $J_\ast$ via left Kan extension.

As in \cite{inverses}, define $IM[k]$ to be the functor $\Tmi \rightarrow \SSets$ given by
\[ F_n \mapsto \Hom_{\Tmi}(F_k, F_n)=(F_k)^n, \] and let $H = I\nerve(F_k)^t$.  The following results continue to hold, replacing $\Tg$ with $\Tmi$.

\begin{lemma} \cite[4.9]{inverses}
In $\LSSets^{\Tmi}_\ast$, $L_2J_*(H)$ is weakly equivalent to $IM[k]$.
\end{lemma}

\begin{prop} \cite[4.10]{inverses}
For any object $X$ in $\LSSets^{\IDeltaop}_{\ast, c}$, $L_1X$ is weakly equivalent to $J^*L_2J_*X$.
\end{prop}

\begin{prop} \cite[4.11, 4.3]{inverses}
The adjoint pair
\[ \xymatrix@1{J_*: \LSSets^{\IDeltaop}_{*,f} \ar@<.5ex>[r] & \SSets^{\Tmi}_* :J^* \ar@<.5ex>[l]} \]
is a Quillen equivalence.
\end{prop}

The more general fixed-object case, as well as the more general case as developed in \cite{inverses2}, can be corrected by considering categories with involution rather than groupoids whenever $\IDeltaop$ is the indexing category.